\input epsf
\epsfverbosetrue

\input amssym.def
\baselineskip=16pt
\hoffset 0truein

\def\bs{\bigskip}
\def\ms{\medskip}
\def\ss{\smallskip}
\def\Z{{\Bbb Z}}
\def\ni{\noindent}
\def\Si{\Sigma}
\def\G{\Gamma}
\def\s{\sigma}
\def\<{\langle}
\def\>{\rangle}
\def\a{\alpha}
\def\b{\beta}
\def\g{\gamma}
\def\d{\delta}
\def\L{\Lambda}
\def\D{\Delta}
\def \={ {\buildrel \cdot \over =}}
\def\qed{{$\vrule height4pt depth0pt width4pt$}}

\centerline{\bf LIFTING REPRESENTATIONS OF $\Z$-GROUPS}\ss

\centerline{Daniel S. Silver and Susan G. Williams} \bs

{\narrower {\bf ABSTRACT:} Let $K$ be the kernel of an
epimorphism $G \to \Z$, where $G$ is a finitely presented
group. If $K$ has infinitely many subgroups of index
$2,3$ or $4$, then it has uncountably many. Moreover, if
$K$ is the commutator subgroup of a classical knot group $G$,
then any homomorphism from $K$ onto the symmetric group $S_2$
(resp. $\Z_3$) lifts to a homomorphism onto $S_3$ (resp.
alternating group $A_4$).}\bs

\ms
\footnote{} {Both authors partially supported by NSF grants
DMS-0071004 and DMS-0304971.}
\footnote{}{2000 {\it Mathematics Subject Classification.}  
Primary 20E07; secondary 57M27.}

\ni {\bf 1. Introduction.}  Let $G$ be a finitely presented
group with infinite abelianization. Given an epimorphism
$\chi: G \to \Z$, we denote its kernel by $K$. Examples of 
special interest arise in knot theory; if $G$ is the group 
$\pi_1(S^3\setminus k)$ of a knot $k\subset S^3$ and $\chi$ is abelianization, then $K$ is 
the commutator subgroup of $G$. 

In general, $K$ need not be finitely generated. Nevertheless,
the Reidemeister-Schreier method [{\bf LS77}] ensures that it
has a group presentation composed of finitely many families
of generators
$a_j, b_j, \ldots, c_j\quad (j\in
\Z)$  and relators  $r_k,
s_k,\ldots, t_k\quad (k \in \Z)$ such that any
relator in a family can be gotten from any other by shifting
all of the indices of its generators by a constant.
(Conversely, any group $K$ with such a presentation arises 
as a kernel $\chi:G  \to \Z$ for some finitely presented $G$.)
Clearly, $K$ admits a nontrivial $\Z$-action by automorphisms.
The action is the restriction to $K$ of conjugation in
$G$ by a preimage $x \in \chi^{-1}(1)$; actions corresponding
to different preimages are related by an inner automorphism of
$K$. For this reason we call $K$  a {\it  finitely
presented $\Z$-group} (cf. [{\bf Ro96}]). 

In [{\bf SW96}] the authors exploited this special structure,
showing that for any finite group $\Si$, the set of
representations ${\rm Hom}(K,\Si)$ has the structure of a
{\it shift of finite type}, a compact
$0$-dimensional dynamical system completely described by a finite directed graph $\G$;
in particular, there is a bijection between ${\rm Hom}(K,
\Si)$ and bi-infinite paths in $\G$. Techniques
of symbolic dynamics can be used to understand ${\rm
Hom}(K,\Si)$. Details are reviewed in \S2. 

Given any group $K$,  its subgroups of index no greater
than $r$ are in finite-to-one correspondence with
representations
$\rho: K \to S_r$, where $S_r$ is the symmetric group on 
$\{1, \ldots, r\}$. The correspondence can be described by 
$\rho \mapsto \{g\in K \mid \rho(g)(1)=1\}$. The preimage
of a subgroup  of index exactly $r$ consists of $(r-1)!$
transitive representations. By a {\it transitive} representation
we mean a representation  $\rho$ such that $\rho(K)$ operates 
transitively on $\{1, \ldots, r\}$. Note that $K$ contains
finitely (resp. countably, uncountably) many subgroups of
index $r$ if and only if ${\rm Hom}(K, S_r)$ contains
finitely (resp. countably, uncountably) many transitive
representations. 

In [{\bf SW99}] we applied techniques of
symbolic dynamics to study ${\rm Hom}(K, S_r)$. Under the
hypothesis  that $K$ {\sl contains an abelian HNN base for} $G$ (see
\S2) we proved that if $K$ contains infinitely many subgroups
of some finite index $r$, then it contains uncountably many.  Our
first result is that this dichotomy continues to hold even if the hypothesis is removed, provided that 
$r<5$.\bs

\ni {\bf  Theorem 3.4. (Dichotomy)} Let $K$ be a finitely presented
$\Z$-group. If $K$ contains infinitely many subgroups of index
$2,3$ or $4$, then $K$ contains uncountably many subgroups of
that index. \bs

The conclusion does not hold if $r>4$ (see
Example 3.4.) As an immediate consequence of the above theorem
and Corollary 1.3 of [{\bf SW99}] we obtain
\bs

\ni {\bf Corollary 3.5.} Let $K$ be a finitely
presented $\Z$-group. If $K$ contains infinitely many subgroups
of index $r=2,3$ or $4$, then $K$ contains uncountably many
subgroups of any index greater than or equal to $r!$. \bs

When $G$ is a knot group,
topology imposes restrictions on the structure of its
commutator subgroup $K$.   \bs

\ni {\bf Theorem 4.3} Assume that $K$ is the commutator
subgroup of a knot group. Then (i) any representation from $K$
onto $S_2$ lifts to a representation onto $S_3$; (ii)
any representation from $K$ onto
$\Z_3$ lifts to a representation onto the alternating group $A_4$. \bs

The conclusions of Theorem 4.3 do not hold for  arbitrary 
$\Z$-groups (see Example 4.4).

The commutator subgroup of a nontrivial fibered knot is free of rank at least 
two. Such a group maps onto any symmetric group (and consequently
contains subgroups of every index). Does such a conclusion hold
for the commutator subgroup of any {\sl nonfibered} knot? We offer a partial answer (see Corollary 3.9).

Much of this paper was inspired by C. Livingston's paper [{\bf
Li95}]. In it he revisits classical results about lifting knot
group representations from the perspective of obstruction
theory. We found that with the aid of symbolic dynamics many of
the techniques extend in a natural way to knot commutator subgroups. From this approach we
obtain new insights into the structure of ${\rm Hom}(K,\Si)$.
  
We are grateful to J. Scott Carter, Alex Suciu and Jim Davis
for helpful discussions. Also, we thank the referees for providing comments and suggestions that improved the 
paper, and Jim Howie for suggesting Example 3.6.
\bs

\ni {\bf 2. Symbolic dynamics and representation shifts.} Assume
that $G$ is a finitely presented group with epimorphism $\chi:
G \to \Z$. Then $G$ can be described as an {\it HNN extension}
$\langle x, B \mid x^{-1}ax = \phi(a),\ \forall a \in U\rangle$.
Here $x \in \chi^{-1}(1),$ 
while $B$ is a finitely generated
subgroup of $K = {\rm ker}\ \chi$.  The map $\phi$ is an
isomorphism between finitely generated subgroups $U,V$ of $B$.
The subgroup $B$ is an {\it HNN base},
$x$ is a {\it stable letter}, $\phi$ is an {\it amalgamating
map}. Details can be found in [{\bf LS77}].
It is possible to choose $B$ so that it contains any
prescribed finite subset of $K$ (see [{\bf Si96}], for
example). 

Conjugation by $x$ induces an automorphism  of $K$. Letting
$B_j = x^{-j}Bx^j, U_j = x^{-j}Ux^j$ and $V_j= x^{-j}Vx^j, j
\in \Z$, we can express $K$ as an infinite free
product in which each subgroup $V_j$ is amalgamated with $U_{j+1}$:
$$K = \langle B_j\mid  x^{-j}\phi(u)x^j= x^{-j-1}ux^{j+1},\ \forall u \in U,\ \forall j \in
\Z\rangle.$$

When $\Si$ is a finite group, the set ${\rm Hom}(K,\Si)$ can be
described by finite directed graph $\G$. The vertex set
consists of all representations $\rho_0: U \to \Si$, a set that
is finite since $U$ is finitely generated. If $\bar \rho_0$ is
a representation from $B$ to $\Si$, then we draw a directed 
edge labeled by $\bar\rho_0$ from the vertex $\rho_0 =
\bar\rho_0\vert_U$ to the vertex $\rho_0' =
\bar\rho_0\vert_V\circ\phi$. Consider any bi-infinite  path in
$\G$ given by an edge sequence
$$\cdots \bar\rho_{-2}\  \bar\rho_{-1}\  \bar\rho_0\  
\bar\rho_1\  \bar\rho_2\  \cdots.$$
The representations from $B_j$ to $\Si$ defined by $y\mapsto
\bar\rho_j(x^jyx^{-j})$ have a unique common extension $\rho:K
\to \Si$. In this way bi-infinite paths of $\G$ correspond to
elements of ${\rm Hom}(K,\Si)$. 

Let $E$ be a finite set (with discrete topology) and give $E^\Z$ the product topology. We define the {\it shift} homeomorphism $\s$ by $(\s y)_j=y_{j+1}$ for $y=(y_j)\in  E^\Z$.  The dynamical system $(E^\Z,\s)$ is called the {\it full shift} on the symbol set $E$.  Any closed $\s$-invariant subset is a {\it subshift}.  In a slight abuse of notation, we will use the same symbol $\s$ for the restriction of $\s$ to a subshift, and also for subshifts on different symbol sets.  A subshift $Y$ is an  $n$-step {\it shift of finite type} if there is a subset $S$ of $E^{n+1}$, the set of {\it allowed $(n+1)$-blocks,} such that Y consists of precisely those $y=(y_j)$ for which $y_j\cdots y_{j+n}\in S$, for all $j$.  (Following the custom in symbolic dynamics and information theory, we write elements of $E^n$ as words  $y_1\cdots y_n$ instead of $n$-tuples.) A finite directed graph $\G$ determines a 1-step shift of finite type $X_\G$, the {\it edge shift} of $\G$:  the symbol set is the edge set $E$ and $ee'$ is an allowed 2-block if the terminal vertex of $e$ is the initial vertex of $e'$.  Thus $X_\G$ is the set of elements of $E^\Z$ that correspond to bi-infinite paths in $\G$.  The graph constructed in the preceding paragraph presents ${\rm Hom}(K,\Si)$ as a shift of finite type.  As we showed in [{\bf SW96}], the product topology on the shift space coincides with the compact-open topology on ${\rm Hom}(K,\Si)$ and the shift map $\s$ is induced by the map $\s_x:{\rm Hom}(K, \Si)  \to {\rm Hom}(K, \Si)$  given by 
$$(\s_x\rho)(a) = \rho(x^{-1}ax),\ \forall a\in K.$$

Two dynamical systems $(Y,\tau)$ and $(Y',\tau')$ are {\it topologically conjugate}, or simply {\it conjugate},  if there is a homeomorphism $h: Y\to Y'$ with $\tau'\circ h=\tau$.  Although the graph presenting ${\rm Hom}(K,\Si)$ depends on the choice of HNN base $B$ and stable letter $x$, the corresponding shift of finite type is uniquely determined up to conjugacy by $G$ and $\chi$.  As in [{\bf SW96}] we call the pair $({\rm Hom}(K, \Si), \s_x)$ 
the {\it representation shift} of $K$ in $\Si$, denoted by $(\Phi_\Si, \s)$ or more
simply $\Phi_\Si$.  An algorithm, based on the Reidemeister-Schreier method, for obtaining a graph presenting $\Phi_\Si$ from a finite presentation of $G$ is given in [{\bf SW96}]. \bs

\ni {\bf Example 2.1.} Consider the Baumslag-Solitar group 
$G= \langle x,a\mid ax = a^2x\rangle$ together with the
epimorphism $\chi:G \to \Z$ mapping $x\mapsto 1$ and $a\mapsto
0$. The group $G$ has an HNN decomposition such that
$B=U$ is the infinite cyclic subgroup $\langle a \rangle$,
$V=\langle a^2\rangle$ and $\phi(a) = a^2$. The kernel $K$ of
$\chi$ has presentation $\langle a_j\mid  a^2_j=a_{j+1}\
\forall j \in \Z\rangle$. 

In this example $\Phi_{\Z_3}$ has exactly 3 elements: the
trivial representation; the representation $\rho$ mapping 
$a_{2j}\mapsto 1$ and each $a_{2j+1}\mapsto 2$; and the
representation $\sigma \rho$ mapping $a_{2j}\mapsto 2$ and each
$a_{2j+1}\mapsto 1$. The graph $\G$ that describes the
representation shift appears in Figure 1. Here we label the vertex
$\rho_0$ by $\rho_0(a)$. \bs


\epsfxsize=2truein
\centerline{\epsfbox{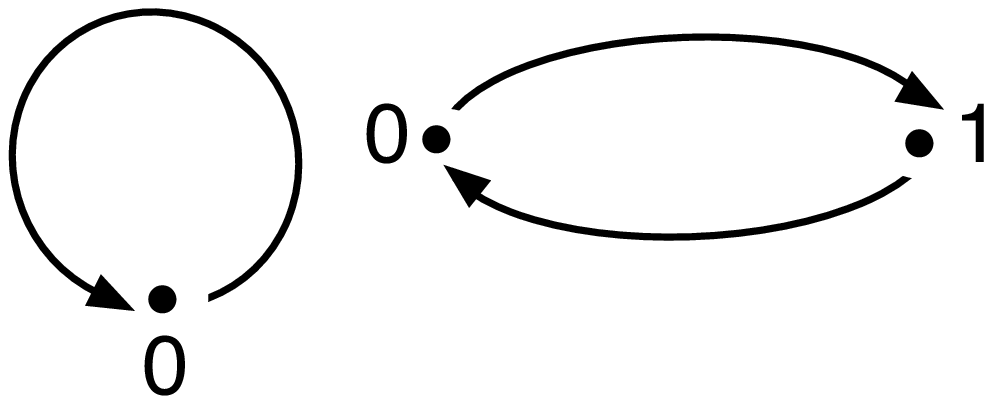}} \bs

\centerline{{\bf Figure 1:} Graph $\G$ describing $\Phi_{\Z_3}$}
\bs

We close this section with some notions from symbolic dynamics that will be needed in the sequel.  We refer the reader to  [{\bf LM95}] for additional background. 

If $Y$ is a subshift of $E^\Z$ and $n$ a positive integer, we may construct a conjugate subshift $Y^{(n)}$ of $(E^n)^\Z$ , the {\it $n$-block presentation} of $Y$, by sending the sequence $(y_j)\in Y$ to the sequence of overlapping $n$-blocks 
$$\ldots(y_j \ldots y_{j+n-1})(y_{j+1}\ldots y_{j+n})(y_{j+2}\ldots y_{j+n+1})\ldots.$$  
If $Y$ is an (at most) $n$-step shift of finite type then $Y^{(n)}$ is a 1-step shift of finite type described by a graph in which the edges are allowed $(n+1)$-blocks of $Y$ and the initial and terminal vertices of an edge are its initial and terminal $n$-blocks.
If we replace the HNN base $B=B_0$ in our construction of a graph presenting $\Phi_\Si$ by the larger base $B^{(n)}=\< B_0\cup B_1\cup\ldots\cup B_{n-1}\>$, we obtain a graph for the $n$-block presentation.

We may, and will, assume that the graph $\G$ of a shift of finite type $X_\G$ has been ``pruned" to remove all vertices and edges that do not lie on a bi-infinite path.  Then $X_\G$ is {\it irreducible} if $\G$ is strongly connected, that is, there is a path from every vertex to every other vertex.  The {\it irreducible components} of $X_\G$ are the finite type subshifts corresponding to the maximal strongly connected subgraphs of $\G$.  

A point $y$ of a subshift has {\it period} $r$ if $\s ^r y=y$.  We will refer to periodic points of the representation shift $\Phi_\Si$ as periodic representations.  Periodic orbits of $X_\G$ correspond to closed paths in $\G$.  If a point $y$ of a shift of finite type has period at most $r$ then it is easily seen that its orbit is represented by a simple closed path in the graph of the $r$-block presentation.  The periodic points of $X_\G$ form a dense subset if and only if $X_\G$ is the union of its irreducible components.  The subshift  $X_\G$ is finite if and only if $\G$ is a disjoint union of cycles, and uncountable if and only if $\G$ has a strongly connected subgraph containing more than one cycle.

A {\it Markov group} is a shift of finite type that is simultaneously a group under an operation that is preserved by the shift map.  If $E$ is a group then the full shift $E^\Z$ is a Markov group under coordinatewise multiplication, and so is any finite type subshift of $E^\Z$ that is also a subgroup.  For example, $\Phi_\Si={\rm Hom}(K, \Si)$ is a Markov group if $\Si$ is abelian.
Theorem 6.3.3 of [{\bf Ki98}] describes the structure of a Markov group: The irreducible component that contains the identity element is topologically conjugate to a full shift (possibly the trivial shift $\{e\}^\Z$), and the entire Markov group is  conjugate to the product of this full shift with a finite Markov group (which may also be trivial).  It follows that a Markov group has dense periodic point set.  It is finite if the full shift factor is trivial, and otherwise every irreducible component is uncountable.

If $K$ is the commutator subgroup of a knot group, then $\Phi_\Si$ must be finite for every abelian $\Si$.  For otherwise, the full shift component would contain a nontrivial fixed point $\rho$.  Then $\rho(a)=\s\rho(a)=\rho(x^{-1}ax)$, or $\rho(a^{-1}x^{-1}ax)=e$, for all $a\in  \rho$.  But the knot group is the normal closure of $x$, and so the elements $a^{-1}x^{-1}ax$ generate $K$
(see [{\bf HK78}], for example.)
\bs


\ni {\bf 3. Lifting representations.} Assume that $E$ is an
extension of
$\Si$ by an abelian group $A$:
$$0\to A \to E\ {\buildrel p \over \to}\
\Si\to 1.$$ For notational convenience, we identify
$A$ with its image
$i(A)$, and regard the latter as a multiplicative group.

The group $E$ acts on $A$ by conjugation. Similarly, there is an
induced action of
$\Si$ on $A$ defined by $a^y := \tilde y a\tilde y^{-1}$, where 
$\tilde y$ is any preimage of $y$; in this way we regard 
$A$ as a $\Si$-module. 
If a homomorphism $\rho: K \to \Si$ lifts to $\tilde \rho: K \to
E$, then $A$ acts by conjugation on the set of liftings. 

The homomorphism $p$ induces a continuous mapping
$p^*:\Phi_E \to \Phi_\Si$. Moreover, $p^*$ is
{\it shift-commuting} in the sense that $p^*\circ \s =
\s\circ p^*$.  If the extension splits then $p^*$ is onto.
\bs

The following proposition is a standard application of 
the cohomology theory of group extensions (see [{\bf Br82}], for
example). \bs

\ni {\bf Proposition 3.1.} Assume that $\rho\in \Phi_\Si$
has a lifting  $\tilde \rho \in \Phi_E$. \ss
(i) The set
of {\it $\rho$-twisted cocycles},
$$ C^1(K; \{A\})=\bigr\{\ \xi:K \to A \mid \xi(xy) =
\xi(x) \xi(y)^{\rho(x)}\ \bigr\},$$
corresponds via $\xi \mapsto \rho_\xi$, where $\rho_\xi(x) =\xi(x) \tilde \rho(x)$, to the complete set
of liftings; \ss
(ii) Liftings $\tilde \rho_1, \tilde \rho_2$ are $A$-conjugate
if and only if $\tilde\rho_2 (\tilde \rho_1)^{-1}$ is a {\it
coboundary}, a map of the form
$x \mapsto a^{\rho(x)}a^{-1}$, for some fixed $a \in A$. 

\bs

\ni {\bf Proof.} 

(i) Given
a derivation $\xi$, define
$\rho_\xi:K \to E$ by $\rho_\xi(x) = \xi(x) \tilde\rho(x)$.
Then
$$\eqalign{\rho_\xi(xy)&= \xi(xy)\tilde\rho(xy) =
\xi(x)\xi(y)^{\rho(x)}\tilde\rho(x)\tilde\rho(y)\cr
&=\xi(x)\tilde\rho(x)\xi(y)\tilde\rho(y)=\rho_\xi(x)\rho_\xi(y).\cr}$$
Hence $\rho_\xi$ is in $\Phi_E$, and it is a lift of $\rho$. 

Conversely, given a lift $\hat\rho$ of $\rho$, define $\xi:K
\to A$ by $\xi(x)=\hat\rho(x)(\tilde\rho(x))^{-1} \in {\rm
Ker}(p) = A$. Then 
$$\eqalign{\xi(xy) &= \hat\rho(xy)(\tilde\rho(xy))^{-1} =
\hat\rho(x)\hat\rho(y)(\tilde\rho(y))^{-1}(\tilde\rho(x))^{-1}\cr
&=\hat\rho(x)(\tilde\rho(x))^{-1}
[\hat\rho(y)(\tilde\rho(y))^{-1}]^{\rho(x)}\cr &= \xi(x)\xi(y)^{\rho(x)}.\cr}$$
Hence $\xi \in C^1(K; \{A\})$ and $\hat \rho = \rho_\xi$.
Clearly, distinct $\xi$ give distinct $\rho_\xi$.

The proof of (ii) is equally routine, and we leave it to the reader.
\qed\bs

\ni {\bf Remark 3.2.} Proposition 3.1 implies that the
$A$-conjugacy classes of liftings correspond to  elements of the
cohomology group $H^1(K;  \{A\})$, with coefficients in $A$
twisted by the action of $K$. 
If $X$ is a complex with $\pi_1X \cong K$, then $H^1(K;  \{A\})$
is isomorphic to $H^1(X, A)$ (see Proposition 2 of
[{\bf Li95}]). 

\bs

\ni {\bf Lemma 3.3.} Let $\rho \in \Phi_\Si$ be a periodic
representation with lifting $\tilde \rho \in \Phi_E$. The
preimage under $p^*$ of the orbit of $\rho$ is a shift of
finite type with dense set of periodic points, that is, a
disjoint union of irreducible shifts of finite type.  The
irreducible components are all finite or all uncountable.   

\bs 

\ni {\bf Proof.} As in \S2, regard $G$ as an HNN extension 
with  finitely generated HNN base $B$ and stable letter 
$x \in \chi^{-1}(1)$. Then
$\Phi_\Si$ is described by a graph with edge set ${\rm Hom}
(B, \Si)$, and $\Phi_E$  by a graph with edge set ${\rm Hom}
(B, E)$. Let $r$ be the least period of $\rho$. Replacing
$B$ with the larger HNN base $B^{(r)}=\< B_0 \cup \cdots \cup
B_{r-1}\>$ if necessary, we can assume that the orbit $O$ of
$\rho$ is represented by a simple cycle in this graph.
Then the set $\tilde O$ of lifts  of $O$ is the finite type
subshift of $\Phi_E$ described by the subgraph consisting of the edges $\eta \in {\rm
Hom}(B,E)$ for which $p\circ \eta$ is an edge in the cycle
presenting $\rho$. 

For $\xi\in C^1(K; \{A\})$ we define $\s\xi: K\to A$ by $\s\xi(y)=\xi(x^{-1}yx)$. It is easy to check that $\s\xi$ is a $\s\rho$-twisted cocycle.  Since $\rho $ has period $r$, we obtain an action of $\tau=\s^r$ on
$C^1(K; \{A\})$. We claim that the pair $(C^1(K; \{A\}), \tau)$
can be viewed as a shift of finite type.  The symbol set is the set $C^1(B^{(r)}; \{A\})$. An element  $\xi\in C^1(K; \{A\})$ can be identified with the sequence of symbols $\xi_j=\tau^j\xi |_{B^{(r)}}$, and $\xi\xi'$ is an an allowed 2-block if $\xi |_{V_{r-1}}\circ\phi=\xi'|_U$.

It is straightforward to check that $C^1(K, \{A\})$ is an
abelian group under the operation $(\xi + \eta)(y)=
\xi(y)\eta(y)$, and that $\tau$ respects this addition.
Hence $(C^1 (K, \{A\}), \tau)$ is a Markov group. By the
Kitchens structure theorem cited in the preceding section, it is a disjoint union of
finitely many shifts of finite type which are all finite
or all uncountable.

By Proposition 3.1, the set ${\cal P}$ of all lifts of $\rho$
is the set of products $\rho_\xi = \xi\tilde\rho$ with 
$\xi \in C^1(K, \{A\})$. Now
$$\s^r(\xi\tilde\rho)(y) = (\xi\tilde\rho)(x^{-r}yx^r) =
\tau\xi(y)\tilde \rho(y),$$
so we can identify the dynamical system $({\cal P}, \s^r)$
with the Markov group $(C^1(K, \{A\}), \tau)$. The shift
of finite type $\tilde O$ is equal to ${\cal P}\cup \s{\cal
P}\cup \cdots \cup \s^{r-1}{\cal P}$. Clearly each irreducible
component of $({\cal P}, \s^r)$ lies in a unique irreducible
component of $(\tilde O, \s)$, and $\tilde O$ is the
disjoint union of these. The components of $\tilde O$ are
finite or uncountable according as the components of ${\cal
P}$ are finite or uncountable. \qed

\bs

\ni {\bf Theorem 3.4.} Assume that $K$ is a finitely presented
$\Z$-group, and $n=2,3$ or $4$. Then $K$ has either finitely
or uncountably many subgroups of index $n$.\bs

\ni {\bf Proof.} It suffices to show that for these $n$, ${\rm
Hom}(K,S_n)=\Phi_{S_n}$ contains either finitely or uncountably
many transitive representations. We may naturally regard $\Phi_{S_2}$ and $\Phi_{S_3}$ as  subshifts of $\Phi_{S_4}$ given by subgraphs of the graph describing  $\Phi_{S_4}$.  

Since $S_2 \cong \Z_2$ is abelian, $\Phi_{S_2}$ is a Markov
group, and hence it is either finite or uncountable. Since all
nontrivial elements are transitive, the theorem holds for $n=2$.

If $\Phi_{S_2}$ is uncountable, then the irreducible component
containing the trivial representation is uncountable. In the
graph describing $\Phi_{S_2}$, the edge corresponding to the trivial representation in ${\rm
Hom}(B,S_2)$ begins and ends at the vertex $v$ corresponding to the trivial representation in ${\rm Hom}(U,S_2)$,  and there must be another path
$p$ from $v$ to itself corresponding to a nontrivial periodic
representation $\rho$. Conjugating $\rho$ with the
transpositions $(23)$ and $(24)$ gives representations $\rho'$
in $\Phi_{S_3}$ and $\rho''$ in $\Phi_{S_4}$ that correspond to
paths $p'$ and $p''$ from $v$ to itself in the graphs describing those representation
shifts. We have $(12)\in \rho(K)$,  and so $(13)\in \rho'(K)$ and $(14)\in \rho''(K).$  Freely concatenating $p$ and $p'$ 
gives uncountably many bi-infinite paths that correspond to transitive representations in $\Phi_{S_3}$; concatenating the paths $p, p'$and $ p''$  gives uncountably many transitive representations in $\Phi_{S_4}$. 

Suppose that $\Phi_{S_2}$ is finite but $\Phi_{S_3}$
is not.  We have a short exact sequence $\Z_3\cong \<(123)\>  \to S_3 \to S_2$ that splits, giving an epimorphism  $\Phi_{S_3}\to \Phi_{S_2}$. There must be
a periodic representation $\rho \in \Phi_{S_2}$ that has
infinitely many lifts in $\Phi_{S_3}$; $\rho$ itself is one of them.  Applying Lemma 3.3, we see that the irreducible component containing $\rho$ in the lift of the orbit of $\rho$ to $\Phi_{S_3}$ is uncountable. Its graph contains a
closed path $p$ corresponding to $\rho$ and another, $\tilde p$, corresponding
to a  lifting $\tilde \rho=\xi\rho$ where $\xi$ is a nontrivial element of $ C^1(K; \{\Z_3\})$. Then either $\tilde \rho(K)$ contains the $3$-cycle
$(123)$, or $\rho(K)$ contains $(12)$ and $\tilde
\rho(K)$ contains $(13)$ or $(23)$. There are uncountably many bi-infinite paths in this component that contain both $p$ and $\tilde p$ and so correspond to transitive representations into $S_3$.  Conjugation by the
transposition $(34)$ fixes $\rho$ but sends $\tilde \rho$ to a
new periodic representation $\bar \rho\in \Phi_{S_4}$, which
must therefore be in the same irreducible component of
$\Phi_{S_4}$ as $\rho$ and $\tilde \rho$. By a similar argument, this component contains uncountably many transitive
representations into $S_4$. 

Finally, suppose that $\Phi_{S_3}$ is finite. Applying Lemma 3.3 to the extension $\Z_2\oplus \Z_2\cong \<(12)(34),(13)(24)\>$ \break $\to S_4\to S_3$ and to each of the periodic orbits of $\Phi_{S_3}$,  we see that
$\Phi_{S_4}$ is a disjoint union of irreducible components, and thus has dense set of periodic points. Each component is finite or uncountable. If $\Phi_{S_4}$ contains infinitely many transitive representations, then so does some
irreducible component.  We can find a closed path $p$ in the graph of this component corresponding to a  periodic transitive representation. Then uncountably many bi-infinite paths in this graph contain $p$ and so correspond to transitive representations.
\qed\bs

\ni {\bf Corollary 3.5.} Let $K$ be a finitely presented
$\Z$-group. If $K$ contains infinitely many subgroups of index
$r = 2,3$ or $4$, then $K$ contains uncountably many subgroups of
any index greater than or equal to $r!$. \bs

\ni {\bf Proof.} Assume that $K$ contains infinitely many
subgroups of index $r = 2,3$ or $4$. Theorem 3.4 implies that there are in fact uncountably many. Consequently, $K$ admits uncountably
many transitive representations into $S_r$. It follows that among
the kernels of the representations there are uncountably many
subgroups of index no greater than $r!$.  Theorem 1.2 of [{\bf
SW99}] implies that $K$ contains uncountably many subgroups
(not necessarily normal) of any index greater than or equal to
$r!$. \qed \bs

It is not difficult to see that the conclusion of Theorem 3.4 holds for arbitrary groups when $r=2$. However, it does not hold for arbitrary groups when $r=3,4$. Nor does it hold 
for $\Z$-groups when $r=5$. The assertion for $r=3$ is established by the following example,
a slight modification of an example generously provided by Jim Howie.\bs

\ni {\bf Example 3.6.}  Every index-3 subgroup $H$ of a group $G$
arises as the preimage $\rho^{-1}(\langle (1 2)\rangle)$, for some transitive representation $\rho: K \to S_3$. Moreover, $H$ is normal if and only if the image of $\rho$ is abelian. 

Let $G$ be the direct sum of countably many copies of $S_3$. Coordinate projections are transitive representations, yielding countably many subgroups of index 3.  We will show that $G$ has no other subgroups $H$ of index 3.

Observe that $H$ cannot be normal. If it were, then $G/H  \cong \Z/3$ would be a quotient of $G/[G,G]$. However, $G/[G,G]$ is a direct sum of countably many copies of ${\Bbb Z}/2$, and any quotient has exponent 2. Hence $H$ arises from a surjective representation 
$\rho: G \to S_3$. We can see that the only such representations are projections composed  with automorphisms of $S_3$. For let $a=(12), b=(123)$, which generate $S_3$. We regard $G$ as generated by $a_i, b_i,\ i\in {\Bbb N}$, 
where $a_i, b_i$ are generators of the $i$th factor of $G$. Since the image of $\rho$ is not abelian, some $b_j$ must be mapped to 
an  element of order 3. The corresponding generator $a_j$ must be sent to an element of order 2. Then no other generator
$a_i,  b_i,\ i\ne j$ can be mapped nontrivially, since such generators commute with $a_j$ and $b_j$ in the direct sum. 

In a similar way, one can show that the direct sum  of countably many copies of $A_4$ contains countably but not uncountably many subgroups of index 4. Details are left to the reader. 

Example 3.6 is notable in another respect. Using a trick of [{\bf St84}] (see page 273), one can see that $G$ is the kernel of homomorphism from a finitely presented group onto ${\Bbb Z}^2$. Hence the conclusion of Theorem 3.4
does not hold for finitely presented ${\Bbb Z}^d$-groups, defined in the obvious manner,  when $d$ is greater than 1. 

We show that the conclusion of Theorem 3.4 does not hold when $r=5$
by modifying Example 3.2 of [{\bf SW99$'$}], an example due to
K.H. Kim and F. Roush.  \bs

\ni {\bf Example 3.7.} The alternating group $A_5$ has presentation
$\<a,b \mid a^2, b^3, (ab)^5\>$. Consider the HNN extension
$G=\langle x, B \mid x^{-1}ax = \phi(a),\ \forall a \in U\rangle$,
where
$$B= \<a, b, a', b' \mid a^2, b^3, (ab)^5, a'^2, b'^3, (a'b')^5, 
[a,a'],\ [b,b'],\ [aa'^{-1}, b']\, [bb'^{-1},
a']\>,$$
a quotient of the free product of two copies of $A_5$;
$U$and $V$ are the subgroups generated by $a,  b$ and  $a', b'$ respectively; and the homomorphism $\phi: U \to V$ maps
$a \mapsto a',  b\mapsto b'$. As  in \S2, let $K$ be the kernel
of the epimorphism $\chi: G \to \Z$ that sends $x \mapsto 1$ and 
$a,b,  a', b' \mapsto 0$.  We will show that the representation shift $\Phi_{S_5}$ is countable.

We construct a graph describing $\Phi_{S_5}$ as in section 2.  A vertex $v\in {\rm Hom}(U,  S_5)$ is determined by $v(a)=\a$, $v(b)=\b$, so we can identify the vertex set with the set of pairs $(\a, \b) \in S_5 \times S_5$ with $\a^2 = \b^3 = (\a\b)^5 = e$.  There is an edge $\bar \rho$ from $(\a, \b)$ to
$(\a', \b')$ if and only if the assignment $\bar \rho(a)=\a$, $\bar \rho(b)=\b$, $\bar \rho(a')=\a'$, $\bar \rho(b')=\b'$ defines an element of ${\rm Hom}(B,  S_5)$.  Since $A_5$ is a simple subgroup
of index 2 in $S_5$, the subgroup generated by  $\a$ and $\b$ is
trivial, cyclic of order $2$ or $A_5$. The second case, in fact,
cannot occur: $\b$ would have to be trivial,
thereby forcing $\a^2 = \a^5 = e$ and so $\a =e$ as well. 

Note that every vertex admits a self-loop (that is, an edge which
begins and ends at the  vertex) and 
an edge from it to the vertex $(e,e)$. As in [{\bf SW99$'$}], we claim that
there are no other edges in $\G$. To see this,  suppose that there
exists an edge from $(e,e)$ to some other vertex $(\a,\b)$. From 
the presentation above, we see that $\a$ and $\b$ must commute, 
and hence $\a = \b=e$.  Now suppose that there exists an edge
from a vertex $(\a,  \b)$ to another $(\g, \d)$, neither of which
is $(e,e)$. Again using the presentation we see that $\a\g^{-1}$
and $\b\d^{-1}$ commute with both $\g$ and  $\d$. Since the center
of $A_5$ is trivial, $\a=\g$ and $\b = \d$ and hence the edge is 
merely a self-loop on $(\a, \b)$. 

It is clear that the graph we have described has countably many bi-infinite paths, and so
the representation shift $\Phi_{S_5}$ is
countable. All of these representations except the trivial one are transitive. Hence $K$ has
countably many subgroups of index $5$ and none of index $2,3$ or $4$. \bs

\ni{\bf Theorem 3.8.} Let $K$ be a finitely presented $\Z$-group. Suppose some representation $\rho$ of $K$ onto $S_2$ has infinitely many lifts to representations into $S_3$.  Then $K$ has uncountably many representations onto $S_n$ for all $n\geq 3$.\bs

\ni {\bf Proof.} The argument is similar to that of the next-to-last paragraph of the proof of Theorem 3.4.  By Lemma 3.3, the component of $\rho$ in the lifting of $\rho$ to $S_3$ is uncountable.  The graph of this component contains closed paths $p$ and $\tilde p$ corresponding to $\rho$ and a nontrivial periodic lifting $\tilde\rho$. The image $\tilde \rho(K)$ must contain $(13)$, $(23)$ or $(123)$, so paths that contain both $p$ and $\tilde p$ correspond to representations onto $S_3$.

Now fix $n \ge 4$. For $4\leq m\leq n$ we can obtain a periodic representation $\rho^{(m)}$ of $K$ onto the subgroup of $S_n$ consisting of all permutations of $\{1,2,m\}$ by conjugating $\tilde\rho$ by the 2-cycle $(3m)$.  Since each of these conjugations fixes $\rho$, they all leave the component of $\rho$ in the graph of $\Phi_{S_n}$ invariant.  Any path in this component that contains the closed paths corresponding to $\rho$, $\tilde\rho$ and each $\rho^{(m)}$ must map onto $S_n$, and there are uncountably many such paths.
\qed\bs

\ni{\bf Corollary 3.9.} Let $K$ be the commutator subgroup of a knot group.  If $K$ has infinitely many representations into $S_3$, then $K$ has uncountably many representations onto $S_n$ for all $n\geq 3$, and hence uncountably many subgroups of every index $n \ge 3$. \bs

\ni{\bf Proof.}  As we noted at the end of Section 2, the Markov groups $\Phi_{\Z_2}$ and $\Phi_{\Z_3}$ are finite in this case.  Hence there are infinitely many representations onto $S_3$, and infinitely many are lifts of a single representation onto $S_2$. \qed\bs

\ni{\bf 4. Additional application to knot groups.}
Assume that $M \cong R^n/A R^m$ is a finitely generated module
over a Noetherian ring
$R$, where the presentation matrix $A$ is an
$n
\times m$ matrix  with entries in $R$. By adjoining zero
columns, we can assume that $m\ge n$. The elementary ideals
$E_i$ of $A$ form a sequence of invariants of $m$. The ideal
$E_i$ is generated by the $(n-i)\times (n-i)$  minors of the matrix
$A$. When $R$ is a factorial domain (for example, $\Z[t,
t^{-1}]$ or $F[t, t^{-1}]$, where $F$ is a field), each
$E_i$ is contained in a unique minimal principal ideal; a
generator, which is well defined up to multiplication by units in $R$, is denoted by $\D_i(t)$. If $R = \Z [t, t^{-1}]$,  then $\D_i$ is a polynomial, 
the $i$th {\it characteristic polynomial}  of 
$M$, and we normalize it so that it has 
the form $c_0 + c_1 t + \cdots + c_d t^d$, where $c_0\ne 0$.

Here we are concerned only with the $0$th characteristic
polynomial.  Note that when $A$ is a
square matrix, $\D_0(t)$ is simply the  determinant of $A$.
Such is the case for any knot: If $X = S^3 \setminus k$ is
a knot complement, and $\tilde X$ is its infinite cyclic
cover, then $H_1(\tilde X; \Z)$ is a finitely generated 
$\Z[t, t^{-1}]$-module with square presentation matrix.
The $0$th characteristic polynomial  is called
the {\it Alexander polynomial} of $k$, denoted here by 
$\D(t)$. The Alexander polynomial $\D(t)$ of any knot has
even degree and satisfies $\D(1) = \pm 1$. 

For any positive integer $r$, we can regard $M$ as a 
finitely generated module over $\Z[s, s^{-1}]$, where 
$s= t^r$. We can obtain a presentation matrix $A(C_r)$ from $A$
by replacing each polynomial entry $f(t)$ by $f(C_r)$, where $C_r$ is the $r\times r$ companion matrix of
the polynomial $t^r -s$: 
$$C_r= \pmatrix{0&0&\cdots&0&s\cr
1&0&\cdots & 0&0\cr
0&1 &\cdots &0&0\cr
\vdots&\vdots &&\vdots&\vdots \cr
0 &0 &\cdots&1&0\cr}$$

\ni {\bf Lemma 4.1.} Assume that $M\cong R^n/ AR^n$ is a
finitely generated $R= \Z[t,t^{-1}]$-module with $0$th
characteristic polynomial
$\D(t) = c_d \prod (t-\a_i)$. Let $s= t^r, r\ge 1$. The $0$th
characteristic polynomial of
$M$ regarded as a $\Z[s, s^{-1}]$-module is $\tilde \D(s) = c_d^r
\prod(s-\a_i^r)$. \bs

\ni {\bf Proof.}  Regard $\Z[s, s^{-1}]$ as a subring of $\Z[t,
t^{-1}]$, which in turn is a subring of ${\bf C}[t, t^{-1}]$.
The matrix $C_r$ is similar over  ${\Bbb C}[t, t^{-1}]$ to 
the diagonal matrix $D={\rm Diag}(t, \zeta t, \ldots,
\zeta^{r-1} t)$, where $\zeta$ is a primitive $r$th root of
unity. Consequently, $\tilde \D(s) = {\rm Det}A(D)$, where as above $A(D)$ denotes the matrix obtained from $A$ by replacing each polynomial entry $f(t)$ by $f(D)$.  By Theorem 1 of [{\bf
KSW99}] the  determinant is equal to
${\rm Det}[\D(D)]$, which is ${\rm Det}[c_d\prod_i (D-\a_i I)]
= c_d^r
\prod_i \prod_j (\zeta^j t - \a_i)$. The last product can be
rewritten as $c_d^r \prod_i \prod_j (t- \zeta ^{-j} \a_i)$,
which is  equal to $c_d^r \prod_i (t^r-\a_i^r)=c_d^r \prod_i
(s-\a_i^r)$. \qed\bs
 
A polynomial $f(t) \in \Z[t, t^{-1}]$ of degree $d$ is {\it
symmetric} (or {\it reciprocal}) if $f(t^{-1})\ \=\  f(t)$. Here ${\buildrel \cdot \over =}$ indicates equality up to multiplication by a unit $\pm t^n$ in $\Z[t, t^{-1}]$. It is well known that 
the Alexander polynomial $\D(t)$ of any knot is a symmetric
polynomial of even degree (see [{\bf BS03}], for example).  \bs

\ni {\bf Corollary 4.2.} Assume the hypotheses of Lemma 4.1.
If, moreover, $\D(t)$ is symmetric, then so is $\tilde\D(s)$. \bs

\ni {\bf Proof.} A polynomial is symmetric if and only if 
its set of zeros (with multiplicities) is sent to itself by
inversion. The desired conclusion follows immediately. \qed \bs

\ni {\bf Theorem 4.3.} Assume that $K$ is the commutator
subgroup of the group of a knot. Then (i) any representation from
$K$ onto $S_2$ lifts to a representation onto $S_3$; (ii)
any representation from $K$ onto
$\Z_3$ lifts to a representation onto $A_4$.\bs

\ni {\bf Proof.} (i) Let  $\rho$ be a representation from $K$
onto $S_2$.
The (unique) epimorphism $S_3 \to S_2$ fits into a short exact
sequence $A=\<(123)\> \to S_3 \to S_2$ that splits, and
hence a lifting $\tilde \rho$ can be found. By Proposition 3.1,
the complete set of liftings $\rho_\xi$ corresponds bijectively to 
the group of twisted cocycles $C^1(K, \{A\})$. It
is not difficult to see $\rho_\xi$ fails to be surjective if
and only if  $\xi$  is a coboundary. We must prove that
$H^1(K, \{A\})$ is  nontrivial. 

Our proof is topological. We construct a 2-complex $X$ with fundamental group
$K$, and invoke Therorem 4 of [{\bf Li95}], an application of Shapiro's Lemma, to see that the
problem of showing that $H^1(K, \{A\})$ is nontrivial is equivalent to proving that the mod-3 first homology group of a certain cover of $X$ has larger rank than then corresponding homology group of $X$. The fact that $K$ is the commutator
subgroup of a knot ensures that the cover satisfies a type of Poincar\'e duality, imposing 
conditions that are sufficient to complete the argument.

Let $X$ be a CW complex with a single vertex and
fundamental group $K$. For any prime $p$,  the   
homology group $H_1(X; \Z_p)$ is a
finitely generated module over the ring $\Z_p[t, t^{-1}]$ of
Laurent polynomials with coefficients in $\Z_p$. A square matrix
presenting the module can be found, and its determinant
$\D_0(H_1(X; \Z_p))$ is the Alexander
polynomial $\D(t)$ of  the knot with coefficients reduced
modulo $p$.  Alternatively, $H_1(X; \Z_p)$ can be viewed as 
a finite-dimensional vector space over $\Z_p$. Its dimension
is equal to the mod p degree of $\D(t)$.  All of the above
statements hold as well using cohomology.

Let $\pi: \tilde X \to X$ be the $2$-fold
cover corresponding to $\rho$. By Theorem 4 of [{\bf Li95}],
$H^1(K, \{A\})\cong H^1(\tilde X; \Z_3)
/\pi^*H^1(X; \Z_3)$.  We must prove that the dimension of 
$H^1(\tilde X; \Z_3)$ exceeds that of $H^1(X;
\Z_3)$. The Universal Coefficient Theorem implies that
the dimension of $H^1(\tilde X; \Z_3)$ (resp. $H^1(X; \Z_3)$) is 
equal to that of $H_1(\tilde X; \Z_3)$ (resp. $H_1(X; \Z_3)$). We will
work with homology. 

Any representation from $K$ to a finite abelian group is 
periodic [{\bf SW99$'$}]. Let $r$ be the period of $\rho$. 
Then $H_1(\tilde X; \Z)$ is a finitely generated
module over $\L= \Z[s, s^{-1}]$, where $s= t^r$. 
In fact $\rho$ induces a homomorphism from the fundamental
group of the $r$-fold cyclic cover $X_r$ of the knot to $S_2$. (Details 
can be found in [{\bf SW99$'$}].) We can
regard $\tilde X$ as an infinite cyclic cover of the induced 2-fold cover of $X_r$.
It follows from Blanchfield duality that $\D_0(H_1(\tilde X; \Z))$ is 
a symmetric polynomial ([{\bf Tu01}], Corollary 14.7). 

We regard $H_1(X; \Z)$ too as a $\L$-module with $0$th characteristic polynomial $\tilde \D(s)$. By Lemma 4.1 and Corollary 4.2, $\tilde \D(s)$ has the same degree as $\D(t)$ and is symmetric.  

The  CW complex $\tilde X$ has a $0$-skeleton $\tilde X^0$ consisting of two vertices. It is convenient to work with the relative homology group $H_1(\tilde X, \tilde X^0; \Z)$; it fits into a short exact sequence $$0\to H_1(\tilde X; \Z) \to H_1(\tilde X, \tilde X^0; \Z) \to \L/(s-1)\to 0,$$ from which it follows that  $$ \D_0(H_1(\tilde X, \tilde X^0;\Z))\  \=\  \D_0(H_1(\tilde X; \Z))(s-1)$$ (see, for example, Lemma 7.2.7 [{\bf Ka96}]).  Consider then  the chain complex $$0\to C_2(\tilde X, \tilde X^0; \Z)\ {\buildrel \partial \over
\to} \ C_1(\tilde X, \tilde X^0; \Z)\to 0$$
for the pair $(\tilde X, \tilde X^0)$.  The boundary $\partial$
can be represented by a matrix of the form
$$ T= \pmatrix{A & B \cr B & A\cr}.$$
Here $A$ and $B$ are square matrices of the same size. The first
half of the columns of $T$ correspond to edges of $\tilde X$ that
are lifts of edges of $X$ beginning at a fixed vertex of $\tilde
X^0$; the remaining columns correspond to edges that are lifts
beginning at the other vertex. 

Row and column operations convert $T$ into 
$$T'= \pmatrix{A-B & B \cr 0 & A+B\cr}.$$ 
The matrix $A+B$ is a relation matrix for $H_1(X; \Z)$, 
and consequently its determinant is $\tilde \D(s)$. 
Since $\D_0(H_1(\tilde X; \Z))$ and $\tilde \D(s)$ are both symmetric polynomials,
so is ${\rm det}(A-B)$. The latter factors as $(s-1)g(s)$ for some (necessarily symmetric) polynomial $g(s)$, since $s-1$ cannot divide $\tilde \D(s)$. 

Clearly ${\rm det}(A-B)$ is congruent modulo 2 to ${\rm det}(A+B)$. It follows
that $g(s)$ has odd degree.  If $g(s)$ (mod 3) is nonzero,
then it must have positive degree. Therefore,
if  $\D_0(H_1(\tilde X; \Z_3)$ is nonzero, then its degree is  larger
than that of $\D_0(H_1(X; \Z_3)$. Equivalently,  if $H_1(\tilde X; \Z_3)$ has
finite dimension, then its dimension is greater than that of $H_1(X; \Z_3)$. \bs

The proof of (ii) follows a similar line of reasoning as (i). 
Let $\rho$ be a representation of
$K$ onto $\Z_3$. The alternating group $A_4$ maps onto $\Z_3$
with kernel $A =\<(12)(34),(13)(24)\>\cong
\Z_2\oplus \Z_2$, and the resulting short exact sequence
splits. Hence a lifting $\tilde \rho: K \to A_4$ can be found.
As in the proof of (i),  Proposition 3.1 implies that
the complete set of liftings $\rho_\xi$ corresponds bijectively to 
the group of twisted cocycles $C^1(K, \{A\})$. Again
$\rho_\xi$ fails to be surjective if and only if  $\xi$  is a
coboundary. We must prove that $H^1(K, \{A\})$ is  nontrivial.

 Let $p: \tilde X \to X$ be the $3$-fold
cover corresponding to $\rho$. By Theorem 4 of [{\bf Li95}], 
$H^1(K, \{A\})\cong H^1(\tilde X; \Z_2)
/p^*H^1(X; \Z_2)$. 
Hence it  suffices
to show that the dimension of  $H_1(\tilde X; \Z_2)$ exceeds 
the dimension of $H_1(X; \Z_2)$. 

Since the unique
vertex of $X$ is covered by 3 vertices in $\tilde X$, a
short exact homology sequence similar to the one above
shows that 
$$\D_0(H^1(\tilde X, \tilde X^0;
\Z_2)) \ \=\ \D_0(H^1(\tilde X; \Z_2)) (s-1)^2.$$

The representation $\rho: K \to \Z_3$ is necessarily periodic,
say of period $r$.  Let $s = t^r$. Then $H_1(\tilde X, \tilde
X^0; \Z_2)$ is a finitely generated $\L = \Z[s,
s^{-1}]$-module. We view $H_1(X, X^0; \Z_2)$ likewise as a
$\L$-module, and again by Blanchfield Duality and Corollary 4.2 its $0$th characteristic polynomial $\tilde \D(s)$ is symmetric of even degree.

As in the proof of (i) we consider  
the chain complex
$$0\to C_2(\tilde X, \tilde X^0; \Z)\ {\buildrel \partial \over
\to} \ C_1(\tilde X, \tilde X^0; \Z)\to 0$$
for the pair $(\tilde X, \tilde X^0)$. The boundary
operator $\partial$ is represented by a matrix of the form
$$T= \pmatrix{A & B &C \cr B & C&A \cr C&A&B\cr},$$
where $A,B$ and $C$ are square matrices of the same size. The
first third of the columns of $T$ correspond to edges of
$\tilde X$ that are lifts of edges of $X$ beginning at a fixed
vertex of $\tilde X^0$; the second third (resp. last third)
correspond to edges that are lifts beginning at the second
(resp. third) vertex. 

The matrix $T$ is similar in $\Z[\zeta][s, s^{-1}]$ to
$$T'= \pmatrix{A+B+C & 0 &0 \cr 0 & A + \zeta B + \zeta^2 C &0
\cr 0&0&A + \zeta^2 B + \zeta C\cr},$$
where $\zeta$ is a primitive 3rd root of unity. Consequently,
$\D_0(H_1(\tilde X, \tilde X^0; \Z))$ is equal to ${\rm
Det}(R')= {\rm Det}(A+B+C) {\rm Det} (A + \zeta B + \zeta^2
C){\rm Det}(A +
\zeta^2 B + \zeta C)$, which we write as 
$\tilde \D(s) F(s) \bar F(s)$. As before $\tilde \D(s)$ is obtained from the Alexander polynomial of the knot,
using Lemma 4.1; $F(s)$ is
a polynomial with coefficients in $\Z[\zeta]$, and $\bar F(s)$
is the polynomial obtained from $F(s)$ by replacing each 
coefficient by its conjugate. 

As in the proof of (i),  $\D_0(H_1(\tilde X; \Z))$ is a symmetric 
polynomial, and hence so is $\D_0(H_1(\tilde X,
\tilde X^0; \Z))$.  Since
$\D_0(H_1(\tilde X, \tilde X^0; \Z)) = \tilde \D(s) F(s) \bar F(s)$
and $\tilde \D(s)$ are symmetric, so is $F(s)\bar F(s)$. By
the Universal Coefficient Theorem, 
$\D_0(H_1(\tilde X, \tilde X^0; \Z_2))$ is equal to
$\tilde \D(s) F(s) \bar F(s)$ with coefficients reduced
modulo 2.

The product $F(s) \bar F(s)$ has integer coefficients. Since the unique extension of the natural projection $\Z \to \Z_3$ to $\Z[\zeta]$ sends $\zeta$ to $1$, the
product $F(s)
\bar F(s)$ is congruent modulo 3 to $(\tilde \D(s))^2$. After suitable
normalization, we can write $F(s) = c_d s^d + \cdots + c_1 s +
c_0$, where each $c_i \in \Z[\zeta]$. The product 
$F(s) \bar F(s)$ is equal to $c_d\bar c_d s^{2d} + \cdots + c_0
\bar c_0.$ If the degree of $F(s)\bar F(s)$ decreases when 
coefficients are reduced modulo 3, then both $c_d \bar c_d$ and
$c_0 \bar c_0$ must be divisible by 3. This implies that $c_d \bar c_d$ and
$c_0 \bar c_0$
vanish under  $\Z[\zeta] \to \Z_3$, and hence both
$c_d$ and $c_0$ also vanish. An induction argument shows that 
the degees of $F(s) \bar F(s)$ and $F(s) \bar F(s)\ ({\rm
mod}\ 3)$ differ by a multiple of 4. But since 
$\tilde \D(s)$ has even degree, the degree of $(\tilde \D(s))^2$ is a
multiple of 4. Hence the degree of $F(s) \bar F(s)$ is a
multiple of 4, and so the degree of $F(s)$ is even. 

Now consider $F(s) \bar F(s)$ with coefficients reduced
modulo 2. As before, if the reduction causes its degree to
decrease, then both $c_d
\bar c_d$ and
$c_0\bar c_0$ must be even.  This implies that both vanish
under the natural projection
$\Z[\zeta] \to \Z_2[\zeta]$. Hence $F(s)
\bar F(s)\ ({\rm mod}\ 2)$
has degree divisible by 4. If it is zero, then $H_1(\tilde X; \Z_2)$
is infinite, and so its dimension is larger than that of $H_1(X;
\Z_2)$. If it is nonzero, then it contains a nontrivial factor
other than
$(s-1)^2$, and  again the dimension of $H_1(\tilde X; \Z_2)$ is
greater than that of
$H_1(X; \Z_2)$. \qed\bs

\ni{\bf Example 4.4.}  It  is easy to construct examples of general $\Z$-groups for which the conclusions of Theorem 4.3 do not hold. For example, consider the 
group $S_2$ presented in the following way.
$$K=\< a_j \mid a_j^2,\ a_{j+1}=a_j,\ \forall j\  \>.$$
Clearly $K$ admits a (unique) homomorphism onto $S_2$ that does not lift to $S_3$.
A similar example can be constructed to show that the second conclusion of Theorem 4.3 does not hold for general finitely presented $\Z$-groups. 

Nontrivial examples can also be constructed. Consider the $\Z$-group:
$$K = \<a_j\mid a_{j+1} = a_j^3\ \forall j\ \>.$$
This group admits a representation $\rho$  onto $S_2$,  mapping each generator
$a_j$ to the nontrivial element. The reader can verify that the matrix $T$ in the proof of Theorem 4.3 is 
$$\pmatrix{s-2& -1 \cr -1 & s-2}.$$
The determinant modulo 3 is equal up to a multiplicative unit to $s-1$, and hence the polynomial $g(s)$ is trivial. Consequently, $\rho$ does not lift onto $S_3$. 

For the second part of Theorem 4.3, consider the group
$$K=\<a_j\mid a_{j+1} = a_j^2\ \forall j\ >$$
of Example 2.1. 
This is the commutator subgroup of a $2$-knot group. It
admits a representation onto $\Z_3$ mapping each $a_{2j}$ to $2$ and each 
$a_{2j+1}$ to $1$. The reader can verify that the matrix $T$ in the proof of Theorem 4.3 is 
$$\pmatrix{1&-1&1&0&0&0 \cr -s&1&0&0&0&1\cr 0&0&1&-1&1&0\cr 0&1&-s&1&0&0\cr 1&0&0&0&1&-1\cr 0&0&0&1&-s&1}.$$
The determinant modulo 2 is equal up to a multiplicative unit to $(s-1)^2$, and hence $F(s)\bar F(s)$ is trivial. Consequently, $\rho$ does not lift onto $A_4$.\bs

\ni {\bf 5. Conclusion.} The obstruction theory used here has provided new insight into the structure of representation shifts. However, questions remain. \bs

\ni {\bf Conjecture 5.1.} Dichotomy holds for commutator subgroups of knot groups. That is, for any finite target group $\Si$, the shift $\Phi_\Si$ is either finite
or uncountable. \bs

A more subtle, dynamical conjecture which would imply Conjecture 5.1  is the following. \bs

\ni {\bf Conjecture 5.2.} For any finite group $\Si$, periodic points are dense in every representation shift $\Phi_\Si$ of a knot commutator subgroup. \bs

A classical result of Perko [{\bf P76}] says that every representation of a knot group onto $S_3$ lifts onto $S_4$. \bs

\ni {\bf Conjecture 5.3.} Let $K$ be the commutator subgroup of a knot $k$.  Every representation of $K$ onto $S_3$ lifts onto $S_4$. \bs

\ni Conjecture 5.3 holds for fibered knots, since in such a case
$K$ is a nonabelian free group. The general conjecture is equivalent to the assertion that for any surjection $\rho: K \to S_3$, the rank of the mod-2 first homology group of
the associated 3-fold dihedral cover $\tilde Y$ of the infinite cyclic cover $Y$ of the knot exceeds that of $Y$.

\bs
\ni {\bf References.} \ms

\ni [{\bf Br82}] K.~S. Brown, Cohomology of Groups, Springer,
New York, 1982. \ss

\ni [{\bf BS03}] G. Burde and H. Zieschang, Knots, 2nd edition, de Gruyter, Berlin, 2003. \ss

\ni [{\bf HK78}] J.C. Hausmann and M. Kervaire, ``Sous-groupes d\'eriv\'es
des groupes de noeuds,'' {\sl L'Enseign. Math.\ \bf 24} (1978), 111--123. \ss

\ni [{\bf Ka96}] A. Kawauchi, A Survey of Knot Theory,
Birkh\"auser Verlag, Basel, 1996. \ss

\ni [{\bf KL99}] P. Kirk and C. Livingston, Twisted Alexander
invariants, Reidemeister torsion, and Casson-Gordon 
invariants, {\sl Topology\ \bf38} (1999), 635--661.\ss

\ni[{\bf Ki98}] B. Kitchens, Symbolic Dynamics: One-sided, Two-sided and Countable State Markov Shifts, Springer-Verlag, Berlin, 1998. \ss

\ni [{\bf KSW99}] I. Kovacs, D.S. Silver and S.G. Williams,
Determinants of commuting-block matrices, {\sl Amer.\ Math.\ Monthly\
\bf106} (1999), 950--952. \ss

\ni [{\bf LM95}] D. Lind, B. Marcus, An Introduction to
Symbolic Dynamics and Coding, Cambridge Univ. Press, Cambridge,
1995. \ss

\ni [{\bf Li95}] C. Livingston, Lifting representations of
knot groups, {\sl J.\ Knot\ Theory\ Ramifications\ \bf4}
(1995), 225--234.\ss

\ni [{\bf LS77}] R.~C. Lyndon and P.~E. Schupp, Combinatorial
Group Theory, Springer, Berlin, 1977. \ss

\ni [{\bf P76}] K.A. Perko, On dihedral covering spaces of knots, {\sl Invent.\ Math.\ \bf 34} (1976), 77--82. \ss

\ni [{\bf Ro96}] D. Robinson, A Course in Group Theory,
Springer, New York, 1996.\ss

\ni [{\bf Si96}] D.~S. Silver, HNN bases and high-dimensional
knots, {\sl Proc.\ Amer.\ Math.\ Soc.\ \bf124} (1996),
1247--1252. \ss

\ni [{\bf SW96}] D.~S. Silver and S.~G. Williams, Augmented
group systems and shifts of finite type, {\sl Israel\ J.\
Math.\ \bf95} (1996), 231--251. \ss

\ni [{\bf SW99}] D.~S. Silver and S.~G. Williams, On groups
with uncountably many subgroups of finite index, {\sl J.\
Pure\ Appl.\ Alg.\ \bf 140} (1999), 75--86. \ss

\ni{[{\bf SW99$'$}]} D.~S. Silver and S.~G. Williams, Knot
invariants from symbolic dynamical systems, {\sl Trans.\
Amer.\ Math.\ Soc.\ \bf 351} (1999), 3243--3265.
\ss
\ni{[{\bf St84}]} R. Strebel, Finitely presented soluble groups, in
Group Theory: Essays for Philip Hall, K.W. Gruenberg and J.E. Roseblade, eds., Academic Press, London, 1984.
\ss

\ni [{\bf Tu01}] V. Turaev, Introduction to Combinatorial
Torsions, Birkh\"auser-Verlag, Basel, 2001. \ss

\ni{[\bf Wa94}] M. Wada, Twisted Alexander polynomial for
finitely presented groups, {\sl Topology\ \bf33} (1994),
241--256. \bs

\ni Department of Mathematics and Statistics, University of
South Alabama, Mobile AL 36688 USA. \ss
\ni e-mail: 
silver@jaguar1.usouthal.edu, swilliam@jaguar1.usouthal.edu

\end